\documentclass{article}%
\usepackage{amsmath}
\usepackage{amsfonts}
\usepackage{amssymb}
\usepackage{graphicx}%
\providecommand{\U}[1]{\protect\rule{.1in}{.1in}}
\newtheorem{theorem}{Theorem}

\newtheorem{corollary}[theorem]{Corollary}

\newtheorem{definition}[theorem]{Definition}
\newtheorem{example}[theorem]{Example}

\newtheorem{remark}[theorem]{Remark}

\begin{document}

\begin{center}
{\LARGE Martingale approximations for random fields}

\bigskip

Magda Peligrad and Na Zhang

\bigskip
\end{center}

Department of Mathematical Sciences, University of Cincinnati, PO Box 210025,
Cincinnati, Oh 45221-0025, USA. \texttt{ }

email: peligrm@ucmail.uc.edu

email: zhangn4@mail.uc.edu

\begin{center}
Abstract
\end{center}

In this paper we provide necessary and sufficient conditions for the mean
square approximation of a random field with an ortho-martingale. The
conditions are formulated in terms of projective criteria. Applications are
given to linear and nonlinear random fields with independent innovations.

\bigskip


\bigskip


\section{\textbf{Introduction}}

\ \ \ \ \ \ A random field consists of multi-indexed random variables
$(X_{u})_{u\in Z^{d}}$. An important class of random fields are
ortho-martingales which were introduced by Cairoli (1969) and have resurfaced
in many recent works. The central limit theorem for stationary
ortho-martingales was recently investigated by Voln\'{y} (2015). It is
remarkable that Voln\'{y} (2015) imposed the ergodicity conditions to only one
direction of the stationary random field. In order to exploit the richness of
the martingale techniques, in this paper we obtain necessary and sufficient
conditions for an ortho-martingale approximation in mean square. These
approximations extend to random fields the corresponding results obtained for
sequences of random variables by Dedecker and Merlev\`{e}de (2002), Zhao and
Woodroofe (2008) and Peligrad (2010). The tools for proving these results
consist of projection decomposition. We present applications of our results to
linear and nonlinear random fields.

We would like to mention several remarkable recent contributions, which
provide interesting sufficient conditions for ortho-martingale approximations,
by Gordin (2009), El Machkouri et al. (2013), Voln\'{y} and Wang (2014), Cuny
et al. (2015), Peligrad and Zhang (2017), and Giraudo (2017). A special type
of ortho-martingale approximation, so called co-boundary decomposition, was
studied by El Machkouri and Giraudo (2017) and Voln\'{y} (2017). Other recent
results involve interesting mixingale-type conditions in Wang and Woodroofe
(2013), and mixing conditions in Bradley and Tone (2017).

Our results could also be formulated in the language of dynamical systems,
leading to new results in this field.

\section{Results}

\ \ \ \ \ \ For the sake of clarity, especially due to the complicated
notation, we shall explain first the results for double indexed random fields
and, at the end, we shall formulate the results for general random fields. No
technical difficulties arise when the double indexed random field is replaced
by a multiple indexed one.

We shall introduce a stationary random field adapted to a stationary
filtration. In order to construct a flexible filtration it is customary to
start with a stationary real valued random field $(\xi_{n,m})_{n,m\in Z}$
defined on a probability space $(\Omega,\mathcal{K},P)$ and to introduce
another stationary random field $(X_{n,m})_{n,m\in Z}$ defined by
\begin{equation}
X_{n,m}=f(\xi_{i,j},i\leq n,j\leq m), \label{defx}%
\end{equation}
where $f$ is a measurable function defined on $R^{Z^{2}}$. Note that $X_{n,m}$
is adapted to the filtration $\mathcal{F}_{n,m}=\sigma(\xi_{i,j},i\leq n,j\leq
m).$ Without restricting the generality we shall define $(\mathbf{\xi
}_{\mathbf{u}})_{\mathbf{u}\in Z^{2}}$ in a canonical way on the probability
space $\Omega$ $=R^{Z^{2}}$, endowed with the $\sigma-$field, $\mathcal{B},$
generated by cylinders. Then, if $\omega=(x_{\mathbf{v}})_{\mathbf{v}\in
Z^{2}}$, we define $\mathbf{\xi}_{\mathbf{u}}^{\prime}(\omega)=x_{\mathbf{u}}%
$. We construct a probability measure $P^{\prime}$ on $\mathcal{B}$ such that
for all $B\in\mathcal{B}$ and any $m$ and $\mathbf{u}_{1},...,\mathbf{u}_{m}$
we have%
\[
P^{\prime}((x_{\mathbf{u}_{1}},...,x_{\mathbf{u}_{m}})\in B)=P((\mathbf{\xi
}_{\mathbf{u}_{1}},...,\mathbf{\xi}_{\mathbf{u}_{m}})\in B).
\]
The new sequence $(\mathbf{\xi}_{\mathbf{u}}^{\prime})_{\mathbf{u}\in Z^{2}}$
is distributed as $(\mathbf{\xi}_{\mathbf{u}})_{\mathbf{u}\in Z^{2}}$ and
re-denoted by $(\mathbf{\xi}_{\mathbf{u}})_{\mathbf{u}\in Z^{2}}$. We shall
also re-denote $P^{\prime}$ as $P.$ Now on $R^{Z^{2}}$ we introduce the
operators%
\[
T^{\mathbf{u}}((x_{\mathbf{v}})_{\mathbf{v}\in Z^{2}})=(x_{\mathbf{v+u}%
})_{\mathbf{v}\in Z^{2}}.
\]
Two of them will play an important role in our paper, namely when
$\mathbf{u=}(1,0)$ and when $\mathbf{u=}(0,1).$ By interpreting the indexes as
notations for the lines and columns of a matrix, we shall call%

\[
T((x_{u,v})_{(u,v)\in Z^{2}})=(x_{u+1,v})_{(u,v)\in Z^{2}}%
\]
the vertical shift and%
\[
S((x_{u,v})_{(u,v)\in Z^{2}})=(x_{u,v+1})_{(u,v)\in Z^{2}}%
\]
the horizontal shift.\ Then define
\begin{equation}
X_{j,k}=f(T^{j}S^{k}(\mathbf{\xi}_{a,b})_{a\leq0,b\leq0}). \label{defXfield}%
\end{equation}

We assume that $X_{0,0}$ is centered and square integrable. We notice that the
variables are adapted to the filtration $(\mathcal{F}_{n,m})_{n,m\in Z}$. To
compensate for the fact that, in the context of random fields, the future and
the past do not have a unique interpretation, we shall consider commuting
filtrations, i.e.
\[
E(E(X|\mathcal{F}_{a,b})|\mathcal{F}_{u,v})=E(X|\mathcal{F}_{a\wedge u,b\wedge
v}).
\]
\ \ \ \ \ This type of filtration is induced, for instance, by an initial
random field $(\xi_{n,m})_{n,m\in Z}$ of independent random variables, or,
more generally can be induced by stationary random fields $(\xi_{n,m})_{n,m\in
Z}$ where only the columns are independent, i.e. $\bar{\eta}_{m}=(\xi
_{n,m})_{n\in Z}$ are independent. This model often appears in statistical
applications when one deals with repeated realizations of a stationary sequence.

It is interesting to point out that commuting filtrations can be described by
the equivalent formulation: for $a\geq u$ we have
\begin{equation}
E(E(X|\mathcal{F}_{a,b})|\mathcal{F}_{u,v})=E(X|\mathcal{F}_{u,b\wedge v}),
\label{pcf}%
\end{equation}
where, as usual, $a\wedge b$ stands for the minimum of $a$ and $b.$ This
follows from this Markovian-type property (see for instance Problem 34.11 in
Billingsley, 1995).

Below we use the notations%
\[
S_{k,j}=\sum\nolimits_{u,v=1}^{k,j}X_{u,v},\ E(X|\mathcal{F}_{a,b}%
)=E_{a,b}(X).
\]
For an integrable random variable $X$, we introduce the projection operators
defined by%
\[
P_{\tilde{0},0}(X)=(E_{0,0}-E_{-1,0})(X)
\]%
\[
P_{0,\tilde{0}}(X)=(E_{0,0}-E_{0,-1})(X).
\]
Note that, by (\ref{pcf}), we have%
\[
{\mathcal{P}}_{{\mathbf{0}}}(X)=P_{\tilde{0},0}P_{0,\tilde{0}}(X)=P_{0,\tilde
{0}}P_{\tilde{0},0}(X).
\]
and by an easy computation and stationarity
\begin{equation}
{\mathcal{P}}_{u,v}(X)=E_{u,v}(X)-E_{u,v-1}(X)-E_{u-1,v}(X)+E_{u-1,v-1}(X).
\label{proj}%
\end{equation}

We shall introduce the definition of an ortho-martingale, which will be
referred to as a martingale with multiple indexes or simply martingale.

\begin{definition}
Let $d$ be a function and define
\begin{equation}
D_{n,m}=d(\xi_{i,j},i\leq n,j\leq m), \label{defD}%
\end{equation}
Assume integrability. We say that $(D_{n,m})_{n,m\in Z}$ is a martingale
differences field if $E_{a,b}(D_{n,m})=0$ for either $a<n$ or $b<m.$
\end{definition}

Set%
\[
M_{k,j}=\sum\nolimits_{u,v=1}^{k,j}D_{u,v}.
\]
In the sequel we shall denote by$\ ||\cdot||$ the norm in $L^{2}.$ By
$\Rightarrow$ we denote the convergence in distribution.

\begin{definition}
We say that a random field $(X_{n,m})_{n,m\in Z}$ defined by (\ref{defx})
admits a martingale approximation if there is a sequence of martingale
differences $(D_{n,m})_{n,m\in Z}$ defined by (\ref{defD})\ such that
\begin{equation}
\lim_{n\wedge m\rightarrow\infty}\frac{1}{nm}||S_{n,m}-M_{n,m}||^{2}%
=0.\label{martapprx}%
\end{equation}

\end{definition}

\begin{theorem}
\label{ThMart Approx}Assume that (\ref{pcf}) holds. The random field
$(X_{n,m})_{n,m\in Z}$ defined by (\ref{defx}) admits a martingale
approximation if and only if%
\begin{equation}
\frac{1}{nm}{\sum\limits_{j=1}^{n}}\sum\limits_{k=1}^{m}||\mathcal{P}%
_{1,1}(S_{j,k})-D_{1,1}||^{2}\rightarrow0\text{ when }n\wedge m\rightarrow
\infty.\label{defDlim2}%
\end{equation}
and both%
\begin{equation}
\frac{1}{nm}||E_{0,m}(S_{n,m})||^{2}\rightarrow0\text{ and }\frac{1}%
{nm}||E_{n,0}(S_{n,m})||^{2}\rightarrow0\text{ when }n\wedge m\rightarrow
\infty.\label{regularity}%
\end{equation}

\end{theorem}

\begin{remark}
\label{Rregularity}Condition (\ref{regularity}) in Theorem \ref{ThMart Approx}
can be replaced by
\begin{equation}
\frac{1}{nm}||S_{n,m}||^{2}\rightarrow||D_{1,1}||^{2}. \label{condvarlim}%
\end{equation}

\end{remark}

\begin{theorem}
\label{ThMartApprox2}Assume that (\ref{pcf}) holds. The random field
$(X_{n,m})_{n,m\in Z}$ defined by (\ref{defx}) admits a martingale
approximation if and only if
\begin{equation}
\frac{1}{nm}{\sum\limits_{j=1}^{n}}\sum\limits_{k=1}^{m}\mathcal{P}%
_{1,1}(S_{j,k})\ \text{converges in }L^{2}\text{ to }D_{1,1}\text{ when
}n\wedge m\rightarrow\infty\label{defDlim}%
\end{equation}
and the condition (\ref{condvarlim}) holds.
\end{theorem}

\begin{corollary}
\label{CorCLT}Assume that the vertical shift $T$ (or horizontal shift $S$) is
ergodic and either the conditions of Theorem \ref{ThMart Approx} or Theorem
\ref{ThMartApprox2} hold. Then
\begin{equation}
\frac{1}{\sqrt{n_{1}n_{2}}}S_{n_{1},n_{2}}\Rightarrow N(0,c^{2})\text{ when
}n_{1}\wedge n_{2}\rightarrow\infty, \label{cltmartf}%
\end{equation}
where $c^{2}=E(D_{0,0})^{2}.$
\end{corollary}

\section{Proofs}

\ \ \ \ \ \ \textbf{Proof of Theorem \ref{ThMart Approx}.}

We start from the following orthogonal representation
\begin{equation}
S_{n,m}={\sum\limits_{i=1}^{n}}\sum\limits_{j=1}^{m}{\mathcal{P}}%
_{i,j}(S_{n,m})+R_{nm},\label{ort dec}%
\end{equation}
with%
\[
R_{n,m}=E_{n,0}(S_{n,m})+E_{0,m}(S_{n,m})-E_{0,0}(S_{n,m}).
\]
Note that for all $1\leq a\leq i-1,$ $1\leq b\leq j-1$ we have ${\mathcal{P}%
}_{i,j}(X_{a,b})=0$; for all $1\leq b\leq j-1$ we have ${\mathcal{P}}%
_{i,j}(X_{i,b})=0$ and for all $1\leq a\leq i-1$, ${\mathcal{P}}_{i,j}%
(X_{a,j})=0$. Whence,
\[
{\mathcal{P}}_{i,j}S_{n,m}={\mathcal{P}}_{i,j}({\sum\limits_{u=i}^{n}}%
\sum\limits_{v=j}^{m}X_{u,v}).
\]
This shows that for any martingale differences sequence defined by
(\ref{defD}), by orthogonality, we obtain%
\begin{align}
||S_{n,m}-M_{n,m}||^{2} &  ={\sum\limits_{i=1}^{n}}\sum\limits_{j=1}%
^{m}||{\mathcal{P}}_{i,j}({\sum\limits_{a=i}^{n}}\sum\limits_{b=j}^{m}%
X_{a,b})-D_{i,j}||^{2}+||R_{n,m}||^{2}\label{decomposition}\\
&  ={\sum\limits_{i=1}^{n}}\sum\limits_{j=1}^{m}||{\mathcal{P}}_{1,1}%
({\sum\limits_{a=1}^{n-i+1}}\sum\limits_{b=1}^{m-j+1}X_{a,b})-D_{1,1}%
||^{2}+||R_{n,m}||^{2}\nonumber\\
&  ={\sum\limits_{i=1}^{n}}\sum\limits_{j=1}^{m}||{\mathcal{P}}_{1,1}%
({S}_{i,j})-D_{1,1}||^{2}+||R_{n,m}||^{2}.\nonumber
\end{align}
A first observation is that we have a martingale approximation if and only if
both (\ref{defDlim2}) is satisfied and $||R_{n,m}||^{2}/nm\rightarrow0$ as
$n\wedge m\rightarrow\infty.$

Computation, involving the fact that the filtration is commuting, shows that%
\begin{equation}
||R_{n,m}||^{2}=||E_{n,0}(S_{n,m})||^{2}+||E_{0,m}(S_{n,m})||^{2}%
-||E_{0,0}(S_{n,m})||^{2},\label{decompose R}%
\end{equation}
and since $||E_{0,0}(S_{n,m})||\leq||E_{0,m}(S_{n,m})||$ we have that
$||R_{n,m}||^{2}/nm\rightarrow0$ as $n\wedge m\rightarrow\infty$ if and only
if (\ref{regularity}) holds. $\square$

\bigskip

\textbf{Proof of Theorem \ref{ThMartApprox2}}

Let us first note that $D_{1,1}$ defined by (\ref{defDlim}) is a martingale
difference. By using the translation operators we then define the sequence of
martingale differences $(D_{u,v})_{u,v\in Z}$ and the sum of martingale
differences $(M_{u,v})_{u,v\in Z}.$ This time we evaluate%
\[
||S_{n,m}-M_{n,m}||^{2}=E(S_{n,m}^{2})+E(M_{n,m}^{2})-2E(S_{n,m}M_{n,m}).
\]
By using the martingale property, stationarity\ and simple algebra we obtain%
\[
E(S_{n,m}M_{n,m})={\sum\limits_{u=1}^{n}}\sum\limits_{v=1}^{m}{\sum
\limits_{i\geq u}^{n}}\sum\limits_{j\geq v}^{m}E(D_{u,v}X_{i,j})={\sum
\limits_{u=1}^{n}}\sum\limits_{v=1}^{m}E(D_{1,1}S_{u,v}).
\]
A simple computation involving the properties of conditional expectation and
the martingale property shows that%
\[
E(D_{1,1}S_{u,v})=E(D_{1,1}\mathcal{P}_{1,1}(S_{u,v})).
\]
By (\ref{defDlim})\ this identity gives that
\[
\lim_{n,m\rightarrow\infty}\frac{1}{nm}E(S_{n,m}M_{n,m})=E(D_{1,1}^{2}).
\]
From the above considerations%
\[
\lim_{n,m\rightarrow\infty}\frac{1}{nm}||S_{n,m}-M_{n,m}||^{2}=\lim
_{n,m\rightarrow\infty}\frac{1}{nm}E(S_{n,m}^{2})-E(D_{1,1}^{2}),
\]
whence the martingale approximation holds by (\ref{condvarlim}).

Let us assume now that we have a martingale approximation. According to
Theorem \ref{ThMart Approx} condition (\ref{defDlim2}) is satisfied. In order
to show that (\ref{defDlim2}) implies (\ref{defDlim}) we apply the
Cauchy-Schwatz inequality twice:%
\begin{align*}
||\frac{1}{nm}{\sum\limits_{i=1}^{n}}\sum\limits_{j=1}^{m}(\mathcal{P}%
_{1,1}(S_{i,j})-D_{1,1})||^{2}  &  \leq\frac{1}{nm^{2}}\sum\limits_{i=1}%
^{n}||\sum\limits_{j=1}^{m}(\mathcal{P}_{1,1}(S_{i,j})-D_{1,1})||^{2}\\
&  \leq\frac{1}{nm}\sum\limits_{i=1}^{n}\sum\limits_{j=1}^{m}||\mathcal{P}%
_{1,1}(S_{i,j})-D_{1,1})||^{2}.
\end{align*}
Also, by the triangular inequality%
\[
|\frac{1}{\sqrt{nm}}||S_{n,m}||-||D_{1,1}||\text{ }|\leq\frac{1}{\sqrt{nm}%
}||S_{n,m}-M_{n,m}||\rightarrow0\text{ as }n\wedge m\rightarrow\infty,
\]
and (\ref{condvarlim}) follows. $\square$

\bigskip

\textbf{Proof of Remark \ref{Rregularity}}

If we have a martingale decomposition, then by Theorem \ref{ThMart Approx} we
have (\ref{defDlim2}) and by Theorem \ref{ThMartApprox2} we have
(\ref{condvarlim}). Now, in the opposite direction, just note that
(\ref{defDlim2}) implies (\ref{defDlim}) and then apply Theorem
\ref{ThMartApprox2}. $\square$

\bigskip

\textbf{Proof of Corollary \ref{CorCLT}. }

This Corollary follows as a combination of Theorem \ref{ThMart Approx} (or
Theorem \ref{ThMartApprox2}) with the main result in Voln\`{y} (2015) via
Theorem 25.4 in Billingsley (1995). $\square$

\section{Multidimensional index sets}

\ \ \ \ \ \ The extensions to random fields indexed by $Z^{d},$ for $d>2,$ are
straightforward following the same lines of proofs as for a two-indexed random
field. By $\mathbf{u\leq n}$ we understand $\mathbf{u=}(u_{1},...,u_{d})$,
$\mathbf{n=}(n_{1},...,n_{d})$ and $1\leq u_{1}\mathbf{\leq}n_{1}$,..., $1\leq
u_{d}\mathbf{\leq}n_{d}\mathbf{.}$ We shall start with a strictly stationary
real valued random field $\mathbf{\xi}=(\xi_{\mathbf{u}})_{\mathbf{u}\in
Z^{d}}$, defined on the canonical probability space $R^{Z^{d}}$ and\ define
the filtrations $\mathcal{F}_{\mathbf{u}}=\sigma(\xi_{\mathbf{j}}%
:\mathbf{j}\leq\mathbf{u})$. We shall assume that the filtration is commuting
if $E_{\mathbf{u}}E_{\mathbf{a}}(X)=E_{\mathbf{u}\wedge\mathbf{a}}(X),$ where
the minimum is taken coordinate-wise. We define
\[
X_{\mathbf{m}}=f((\xi_{\mathbf{j}})_{\mathbf{j}\leq\mathbf{m}})\text{ and set
}S_{\mathbf{k}}=\sum\nolimits_{\mathbf{u}=\mathbf{1}}^{\mathbf{k}%
}X_{\mathbf{u}}.
\]
We also define $T_{i}$ the coordinate-wise translations and then
$X_{\mathbf{k}}=f(T_{1}^{k_{1}}\circ...\circ T_{d}^{k_{d}}(\xi_{\mathbf{u}%
})_{\mathbf{u}\leq\mathbf{0}}).$ Let $d$ be a function and define%
\begin{equation}
D_{\mathbf{m}}=d((\xi_{\mathbf{j}})_{\mathbf{j}\leq\mathbf{m}})\text{ and set
}M_{\mathbf{k}}=\sum\nolimits_{\mathbf{u}=\mathbf{1}}^{\mathbf{k}%
}D_{\mathbf{u}}.\label{def Dg}%
\end{equation}
Assume integrability. We say that $(D_{\mathbf{m}})_{\mathbf{m}\in Z^{d}}$ is
a martingale differences field if $E_{\mathbf{a}}(D_{\mathbf{m}})=0$ is at
least one coordinate of $\mathbf{a}$ is strictly smaller than the
corresponding coordinate of $\mathbf{m.}$ We have to introduce the
$d$-dimensional projection operator. By using the fact that the filtration is
commuting, it is convenient to define%

\[
\mathcal{P}_{\mathbf{0}}(X)=P_{1}\circ P_{2}\circ...\circ P_{d}(X),
\]
where%
\[
P_{j}(Y)=E(Y|\mathcal{F}_{0}^{(j)})-E(Y|\mathcal{F}_{-1}^{(j)}).
\]
Above, we used the notation: $\mathcal{F}_{0}^{(j)}=\mathcal{F}_{\mathbf{0}}$,
and $\mathcal{F}_{-1}^{(j)}=\mathcal{F}_{\mathbf{u}}$, where $\mathbf{u}$ has
all the coordinates $0$ with the exception of the $j$-th coordinate, which is
$-1$. For instance when $d=3,$ $P_{2}(Y)=E(Y|\mathcal{F}_{0,0,0}%
)-E(Y|\mathcal{F}_{0,-1,0}).$

We say that a random field $(X_{\mathbf{n}})_{\mathbf{n}\in Z^{d}}$ admits a
martingale approximation if there is a sequence of martingale differences
$(D_{\mathbf{m}})_{\mathbf{m}\in Z^{d}}$ such that
\begin{equation}
\frac{1}{|\mathbf{n|}}||S_{\mathbf{n}}-M_{\mathbf{n}}||^{2}\rightarrow0\text{
when }\min_{1\leq i\leq d}n_{i}\rightarrow\infty.\label{martapprxd}%
\end{equation}
where $|\mathbf{n|=}n_{1}...n_{k}.$

Let us introduce the following regularity condition%
\begin{equation}
\frac{1}{|\mathbf{n|}}||S_{\mathbf{n}}||^{2}\rightarrow E(D_{\mathbf{1}}%
^{2})\text{ when }\min_{1\leq i\leq d}n_{i}\rightarrow\infty.\label{varlimg}%
\end{equation}

\begin{theorem}
\label{ThMart Approx3}Assume that the filtration is commuting. The following
statements are equivalent: \newline(a) The random field $(X_{\mathbf{n}%
})_{\mathbf{n}\in Z^{d}}$ admits a martingale approximation. \ \newline(b) The
random field satisfies (\ref{varlimg}) and%
\begin{equation}
\frac{1}{|\mathbf{n|}}\sum\nolimits_{\mathbf{j}\geq\mathbf{1}}^{\mathbf{n}%
}||\mathcal{P}_{\mathbf{1}}(S_{\mathbf{j}})-D_{\mathbf{1}}||^{2}%
\rightarrow0\text{ when }\min_{1\leq i\leq d}n_{i}\rightarrow\infty
.\label{defDlimg}%
\end{equation}
\newline(c) The random field satisfies (\ref{defDlimg}) and for all $j$,
$1\leq j\leq d$ we have
\[
\frac{1}{|\mathbf{n|}}||E_{\mathbf{n}_{j}}(S_{\mathbf{n}})||^{2}%
\rightarrow0\text{ when }\min_{1\leq i\leq d}n_{i}\rightarrow\infty.
\]
where and $\mathbf{n}_{j}\mathbf{\in}Z^{d}$ has the $j-th$ coordinate $0$ and
the other coordinates equal to the coordinates of $\mathbf{n}$.\newline(d) The
random field satisfies (\ref{varlimg}) and
\begin{equation}
\frac{1}{|\mathbf{n|}}{\sum\limits_{\mathbf{j}=\mathbf{1}}^{\mathbf{n}}%
}P_{\mathbf{1}}(S_{\mathbf{j}})\ \text{converges in }L^{2}\text{ to
}D_{\mathbf{1}}\text{ when }\min_{1\leq i\leq d}n_{i}\rightarrow
\infty.\label{defDlimg2}%
\end{equation}

\end{theorem}

\begin{corollary}
Assume that one of the shifts $(T_{i})_{1\leq i\leq d}$ is ergodic and either
one of the conditions of Theorem \ref{ThMart Approx3} holds. Then
\[
\frac{1}{\sqrt{|\mathbf{n|}}}S_{\mathbf{n}}\Rightarrow N(0,c^{2})\text{ when
}\min_{1\leq i\leq d}n_{i}\rightarrow\infty,
\]
where $c^{2}=||D_{\mathbf{0}}||^{2}.$
\end{corollary}

\section{Examples}

\ \ \ \ \ \ Let us apply these results to linear and nonlinear random fields
with independent innovations.

\begin{example}
\label{exinear}(Linear field) Let $(\xi_{\mathbf{n}})_{\mathbf{n}\in Z^{d}}$
be a random field of independent, identically distributed random variables
which are centered and have finite second moment. Define%
\[
X_{\mathbf{k}}=\sum_{\mathbf{j}\geq\mathbf{0}}a_{\mathbf{j}}\xi_{\mathbf{k}%
-\mathbf{j}}.
\]
Assume that $\sum_{\mathbf{j}\geq\mathbf{0}}a_{\mathbf{j}}^{2}<\infty$ and
denote $b_{\mathbf{j}}=\sum_{\mathbf{k}=\mathbf{1}}^{\mathbf{j}}a_{\mathbf{k}%
}.$ Also assume that%
\begin{equation}
\frac{1}{|\mathbf{n|}}\sum_{\mathbf{j}=\mathbf{1}}^{\mathbf{n}}b_{\mathbf{j}%
}\rightarrow c\text{ when }\min_{1\leq i\leq d}n_{i}\rightarrow\infty
\label{linear}%
\end{equation}
and%
\[
\frac{E(S_{\mathbf{n}}^{2})}{|\mathbf{n|}}\rightarrow c^{2}\sigma^{2}\text{
when }\min_{1\leq i\leq d}n_{i}\rightarrow\infty.
\]
Then the martingale approximation\ holds.
\end{example}

\textbf{Proof of Example \ref{exinear}}. The result follows by simple
computations and by applying Theorem \ref{ThMart Approx3} (d). $\square$

\bigskip

\begin{example}
\label{Volterra}(Volterra field) Let $(\xi_{\mathbf{n}})_{\mathbf{n}\in Z^{d}%
}$ be a random field of independent random variables identically distributed
centered and with finite second moment. Define%
\[
X_{\mathbf{k}}=\sum_{(\mathbf{u},\mathbf{v)}\geq(\mathbf{0},\mathbf{0}%
)}a_{\mathbf{u},\mathbf{v}}\xi_{\mathbf{k-u}}\xi_{\mathbf{k-v}},
\]
where $a_{\mathbf{u},\mathbf{v}}$ are real coefficients with $a_{\mathbf{u}%
,\mathbf{u}}=0$ and $\sum_{\mathbf{u,v}\geq\mathbf{0}}a_{\mathbf{u,v}}%
^{2}<\infty$ and assume
\[
\frac{E(S_{\mathbf{n}}^{2})}{|\mathbf{n|}}\rightarrow c^{2}\sigma^{2}\text{
when }\min_{1\leq i\leq d}n_{i}\rightarrow\infty.
\]
Denote%
\[
c_{\mathbf{n,w}}=\frac{1}{|\mathbf{n|}}\sum_{\mathbf{j}=\mathbf{1}%
}^{\mathbf{n}}\sum_{\mathbf{k}=\mathbf{1}}^{\mathbf{j}}(a_{\mathbf{k}%
,\mathbf{k-w}}+a_{\mathbf{k-w},\mathbf{k}})
\]
and assume that
\[
\lim_{\mathbf{n>m}\rightarrow\infty}\sum_{\mathbf{w}\leq\mathbf{0}%
}(c_{\mathbf{n,w}}-c_{\mathbf{m,w}})^{2}=0.\
\]
Then the martingale approximation holds.
\end{example}

\textbf{Proof of Example \ref{Volterra}}. We have
\[
\mathcal{P}_{\mathbf{0}}(X_{\mathbf{k}})=\sum_{(\mathbf{u},\mathbf{v)}%
\geq(\mathbf{0},\mathbf{0})}a_{\mathbf{u},\mathbf{v}}\mathcal{P}_{\mathbf{0}%
}(\xi_{\mathbf{k-u}}\xi_{\mathbf{k-v}}).\
\]
Note that $\mathcal{P}_{\mathbf{0}}(\xi_{\mathbf{k-u}}\xi_{\mathbf{k-v}}%
)\neq0$ if and only if $\mathbf{k-u=0}$ and $\mathbf{k-v=w}$ with
$\mathbf{w}\leq\mathbf{0}$ or $\mathbf{k-v=0}$ and $\mathbf{k-u=t\leq0.}$ Therefore

\textbf{ }%
\begin{align*}
\mathcal{P}_{\mathbf{0}}(X_{\mathbf{k}})  &  =\sum_{\mathbf{w}\leq\mathbf{0}%
}a_{\mathbf{k},\mathbf{k-w}}\xi_{\mathbf{0}}\xi_{\mathbf{w}}+\sum
_{\mathbf{t}\leq\mathbf{0}}a_{\mathbf{k-t},\mathbf{k}}\xi_{\mathbf{0}}%
\xi_{\mathbf{t}}\ \\
&  =\sum_{\mathbf{w}\leq\mathbf{0}}(a_{\mathbf{k},\mathbf{k-w}}%
+a_{\mathbf{k-w},\mathbf{k}})\xi_{\mathbf{0}}\xi_{\mathbf{w}}.
\end{align*}
So%
\[
\mathcal{P}_{\mathbf{0}}(S_{\mathbf{j}})=\sum_{\mathbf{w}\leq\mathbf{0}}%
\sum_{\mathbf{k}=\mathbf{1}}^{\mathbf{j}}(a_{\mathbf{k},\mathbf{k-w}%
}+a_{\mathbf{k-w},\mathbf{k}})\xi_{\mathbf{0}}\xi_{\mathbf{w}}\
\]
and%
\[
\frac{1}{|\mathbf{n|}}\sum_{\mathbf{j}=\mathbf{1}}^{\mathbf{n}}\mathcal{P}%
_{\mathbf{0}}(S_{\mathbf{j}})=\frac{1}{|\mathbf{n|}}\sum_{\mathbf{w}%
\leq\mathbf{0}}\sum_{\mathbf{j}=\mathbf{1}}^{\mathbf{n}}\sum_{\mathbf{k}%
=\mathbf{1}}^{\mathbf{j}}(a_{\mathbf{k},\mathbf{k-w}}+a_{\mathbf{k-w}%
,\mathbf{k}})\xi_{\mathbf{0}}\xi_{\mathbf{w}}\ .
\]
By independence and because $a_{\mathbf{k},\mathbf{k}}=0,$ this expression
converges in $L^{2}$ if and only if%
\[
\frac{1}{|\mathbf{n|}}\sum_{\mathbf{w}\leq\mathbf{0}}\sum_{\mathbf{j}%
=\mathbf{1}}^{\mathbf{n}}\sum_{\mathbf{k}=\mathbf{1}}^{\mathbf{j}%
}(a_{\mathbf{k},\mathbf{k-w}}+a_{\mathbf{k-w},\mathbf{k}})\xi_{\mathbf{0}}%
\xi_{\mathbf{w}}\
\]
converges in $L^{2}.$ By independence, and with the notation%
\[
c_{\mathbf{n,w}}=\frac{1}{|\mathbf{n|}}\sum_{\mathbf{j}=\mathbf{1}%
}^{\mathbf{n}}\sum_{\mathbf{k}=\mathbf{1}}^{\mathbf{j}}(a_{\mathbf{k}%
,\mathbf{k-w}}+a_{\mathbf{k-w},\mathbf{k}}),
\]
this convergence happens if
\[
\lim_{\mathbf{n>m}\rightarrow\mathbf{\infty}}\sum_{\mathbf{w}\leq\mathbf{0}%
}(c_{\mathbf{n,w}}-c_{\mathbf{m,w}})^{2}=0.\
\]
It remains to apply Theorem \ref{ThMart Approx3} (d).

\bigskip

\textbf{Acknowledgements.} This research was supported in part by the
NSF\ grant DMS-1512936 and the Taft Research Center at the University of Cincinnati.

\end{document}